\documentclass[12pt,a4paper]{article}

\usepackage[T1]{fontenc}
\usepackage{lmodern}
\usepackage[english]{babel}
\usepackage{amsmath}
\usepackage{amssymb}
\usepackage{booktabs}
\usepackage{hyperref}
\usepackage{microtype}
\usepackage{multirow}
\usepackage{subcaption}

\newtheorem{theorem}{Theorem}
\newtheorem{conjecture}{Conjecture}

\makeatletter
\newcommand{\STS@text}[1]{{\rm STS}$(#1)$}
\newcommand{\STS@unstarred}[1]{\ifmmode{\text{\mbox{\STS@text{#1}}}}\else{\mbox{\STS@text{#1}}}\fi}
\newcommand{\STS@starred}[1]{\ifmmode{\text{\mbox{sub-\STS@text{#1}}}}\else{\mbox{sub-\STS@text{#1}}}\fi}
\newcommand{\STSS@unstarred}[1]{\ifmmode{\text{\mbox{\STS@text{#1}s}}}\else{\mbox{\STS@text{#1}s}}\fi}
\newcommand{\STSS@starred}[1]{\ifmmode{\text{\mbox{sub-\STS@text{#1}s}}}\else{\mbox{sub-\STS@text{#1}s}}\fi}
\newcommand{\STS}{\@ifstar{\STS@starred}{\STS@unstarred}}
\newcommand{\STSS}{\@ifstar{\STSS@starred}{\STSS@unstarred}}
\makeatother

\begin{document}

\title{Steiner Triple Systems of Order 21 with Subsystems}

\author{
Daniel Heinlein\thanks{Supported by the Academy of Finland, Grant 331044.}\phantom{ } and Patric R. J. \"Osterg\aa rd\\
Department of Communications and Networking\\
Aalto University School of Electrical Engineering\\
P.O.\ Box 15400, 00076 Aalto, Finland\\
\tt \{daniel.heinlein,patric.ostergard\}@aalto.fi
}

\date{}

\maketitle

\begin{abstract}
  The smallest open case for classifying Steiner triple systems is
  order 21. A Steiner triple system of order 21, an \STS{21}, can have
  subsystems of orders 7 and 9, and it is known that there are
  12,661,527,336 isomorphism classes of \STSS{21} with \STSS*{9}.
  Here, the classification of \STSS{21} with subsystems is completed by
  settling the case of \STSS{21} with \STSS*{7}.
  There are
  116,635,963,205,551 isomorphism classes of such systems. An estimation
  of the number of isomorphism classes of \STSS{21} is given.
  \end{abstract}



\section{Introduction}

A \emph{Steiner triple system} (STS) is a pair $(V,\mathcal{B})$,
where $V$ is a set of \emph{points} and $\mathcal{B}$ is a
set of \mbox{3-subsets} of points, called \emph{blocks}, such that
every \mbox{2-subset} of points occurs in exactly one block. The size of
the point set is the \emph{order} of the STS, and an
STS of order $v$ is denoted by \STS{v}.
It is well known that an \STS{v} exists iff
\begin{equation}
\label{eq:sts}
v \equiv 1\text{ or }3\!\!\!\pmod{6}.
\end{equation}
For more information about Steiner triple systems, see~\cite{C,CR}.

An \STS{v} is said to be \emph{isomorphic} to another \STS{v} if
there exists a bijection between the point sets that maps blocks
onto blocks; such a bijection is called an \emph{isomorphism}.
An isomorphism of a Steiner triple system onto itself
is an \emph{automorphism} of the STS.
The automorphisms of an STS form a group under composition, the
\emph{automorphism group} of the Steiner triple system.

Classification of combinatorial designs is about finding a transversal
of the isomorphism classes~\cite{KO2}. The Steiner triple systems
have been classified up to order 19, and the numbers of isomorphism classes
are 1, 1, 1, 2, 80, and 11,084,874,829 for orders 3, 7, 9, 13, 15, and
19, respectively. A classification of the \STSS{19} was published in
2004 with a remark that the algorithm used would require hundreds of thousands
of CPU years to classify the \STSS{21}~\cite{KO1}.
As this seems to be currently out of reach, one will have to focus on subclasses of \STSS{21}.
Indeed, \STSS{21} of various types have been considered in this context,
including \STSS{21} with a nontrivial automorphism group~\cite{K1}
(with earlier work in~\cite{CCIL,KT,LM,MPR,MR1,T1,T2}, also
considering other properties), anti-Pasch \STSS{21}~\cite{KO4},
and resolutions of \STSS{21}---that is, Kirkman
triple systems---with subsystems~\cite{KO5}.

A necessary condition for an \STS{v} to have a
nontrivial ($w>3$) and proper ($w<v$) subsystem of order $w$,
i.e., a \STS*{w}, is that
$v \geq 2w+1$; see~\cite[Lemma~6.1]{CR}. Classification of Steiner
triple systems with \STSS*{7} has been carried out for orders
15 and 19---see~\cite[Table~1.29]{MR2} and~\cite{KOTZ}, respectively---and
for those with \STSS*{9} for order 19---see~\cite{SS}.

The only possible nontrivial proper subsystems of \STSS{21} are
\STSS{7} and \STSS{9}. The \STSS{21} with \STSS*{9} are classified
in~\cite{KOP}; there are
12,661,527,336 isomorphism classes of such designs. For \STSS{21} with
\STSS*{7}, the special case of Wilson-type systems is handled in~\cite{KOTZ}.
\emph{Wilson-type} \STSS{21} contain three \STSS*{7} on disjoint
point sets. In the current paper the classification problem for \STSS{21}
with subsystems is settled by completing the case of \STSS*{7}.

\begin{theorem}
There are 116,635,963,205,551 isomorphism classes of \STSS{21}
containing at least one \STS*{7}.
\end{theorem}

The paper is organized as follows.
An algorithm for classifying \STSS{21} with \STSS*{7} is described
in Section~\ref{sect:class}, and
the results are listed in Section~\ref{sect:res}.
The number of isomorphism classes of \STSS{21}
with \STSS*{7} is used in
Section~\ref{sect:estimate} to
get an estimation of the total number of isomorphism classes of
\STSS{21}.

\section{Classification}\label{sect:class}

In this section, we present a classification algorithm for \STSS{21} containing
\STSS*{7}.
To facilitate reading, we give necessary definitions in Section~\ref{sec:defs}.
The general approach is outlined in
Section~\ref{sect:general}, details about subtasks are given
in Section~\ref{sect:detail}, and some computational issues
are considered in Section~\ref{sect:sub}.

\subsection{Definitions}\label{sec:defs}

A $(v_r,b_k)$ configuration is an incidence structure with $v$ points and $b$ blocks, such that each block contains $k$ points, each point occurs in $r$ blocks, and two different blocks intersect in at most one point.
If $v=b$ and $k=r$, these are simply called $v_k$ configurations.
The definitions of isomorphism and automorphism of configurations are analogous to those for
Steiner triple systems.

A \emph{\mbox{1-factor}} in a graph, also called a perfect matching, is a \mbox{1-regular} spanning subgraph and a \emph{\mbox{1-factorization}} is a partition of the edges of the graph into \mbox{1-factors}.
A \mbox{1-factorization} of
a graph $G = (V,E)$ is isomorphic to a \mbox{1-factorization} of a graph $G' = (V',E')$
if there is a bijection from $V$ to $V'$ that maps the
\mbox{1-factors} of the \mbox{1-factorization} of $G$ onto the \mbox{1-factors} of
the \mbox{1-factorization} of $G'$.

\subsection{General Approach}\label{sect:general}

On a general level, the current approach follows~\cite{KOP}, in which all of Theorem~\ref{theo:splitting} except the last statement already appeared.

\begin{theorem}[\cite{KOP}]\label{theo:splitting}
Let $(V,\mathcal{B})$ be an \STS{v} that has a \STS*{w} $(W,\mathcal{B}')$.
Then
\begin{enumerate}
\item $\mathcal{B} = \mathcal{B}' \cup \mathcal{F} \cup \mathcal{D}$ where
$\mathcal{F}$ and $\mathcal{D}$ are the sets of blocks that intersect $W$ in 1 and 0 points, respectively,
\item $\mathcal{F} = \bigcup_{p \in W} \mathcal{B}_p$ where $\mathcal{B}_p$ is the set of blocks in $\mathcal{F}$ that contain $p \in W$,
\item $\mathcal{B}'_p = \{B \setminus \{p\} : B \in \mathcal{B}_p\}$ is a \mbox{1-factor} of a graph $G$ with vertices $V \setminus W$ and edges $\bigcup_{p \in W} \mathcal{B}'_p$,
\item $\{\mathcal{B}'_p : p \in W\}$ is a \mbox{1-factorization} of $G$,
\item $G$ is $w$-regular and its complement $\overline{G}$ is $(v-2w-1)$-regular, and
\item $\overline{G}$ can be decomposed into a set of edge-disjoint \mbox{3-cycles}---$\mathcal{D}$ being one possible set---which forms a $$({(v-w)}_{(v-2w-1)/2},{((v-w)(v-2w-1)/6)}_{3})$$ configuration.
\end{enumerate}
\end{theorem}

Using this theorem, any STS containing a sub-STS is decomposable into $\mathcal{B}' \cup \mathcal{F} \cup \mathcal{D}$.
For the task of classifying all \STSS{v} containing some \STS*{w}, one has now two starting points:
either a classification of the \mbox{1-factorizations} underlying $\mathcal{F}$ or a classification of the configurations corresponding to $\mathcal{D}$.
Then, in both cases, one needs to combine this with a classification of $\mathcal{B}'$ to create an STS in all possible ways, taking symmetry into account.

The next sections illustrate the details for $v=21$ and $w=7$; the general setting is also depicted in~\cite{KOP}.

\subsection{Application to STS(21) containing sub-STS(7)s}

Let $(V,\mathcal{B})$ be an \STS{21} that has a
\STS*{7} $(W,\mathcal{B}')$. Clearly $W \subseteq V$
and $\mathcal{B}' \subseteq \mathcal{B}$. The blocks in
$\mathcal{B} \setminus \mathcal{B'}$ intersect $W$ in
either 0 or 1 points, and those two sets of blocks are
denoted by $\mathcal{D}$ and $\mathcal{F}$, respectively.

Fix a point $p \in W$ and let $\mathcal{B}_p$ be the
set of blocks in $\mathcal{F}$ that contain $p$. Further
let

\begin{equation}
\label{eq:1}
\mathcal{B}'_p = \{B \setminus \{p\} : B \in \mathcal{B}_p\}.
\end{equation}

As a pair of points with one point in $W$ and
the other in $V \setminus W$ must occur in exactly
one block of $\mathcal{F}$, the sets in $\mathcal{B}'_p$
partition $V \setminus W$. The sets in
$\mathcal{B}'_p$ have size 2, and we may view them as
edges in a graph with vertex set $V \setminus W$.
The sets in $\mathcal{B}'_p$ form a \mbox{1-factor} of
that graph. With 7 possible values of $p$,
we have 7 disjoint \mbox{1-factors} of a
\mbox{7-regular} graph of order 14.

To complete the Steiner triple system,
given a \mbox{7-regular} graph $G$ of order 14, one may
find all \mbox{1-factorizations} of $G$ and in the complement $\overline{G}$
find all decompositions into \mbox{3-cycles} (which is the graph
analogy of finding sets of triples that cover all unordered pairs)
and combine these in all possible ways. Doing this for all
possible choices of $G$ gives all ways of extending the
initial \STS{7}.
Finally, isomorph rejection needs to be carried out during
the process of combining parts.
Specific details about using this approach in the current
work---where the order $\mathcal{D} \rightarrow \mathcal{F} \rightarrow
\mathcal{B}'$ for constructing the blocks $\mathcal{B}$ is actually
used---are presented in Section~\ref{sect:detail}.
See also~\cite{KOP}.

There are 21,609,301 
isomorphism classes of \mbox{7-regular} graphs of order 14~\cite{M2};
see also~\cite[Table~4.25]{R}.
Only a small number of graphs $G$ have the property that the complement
$\overline{G}$ can be decomposed into \mbox{3-cycles}
as described in the last statement in Theorem~\ref{theo:splitting}.
Indeed, the required $14_3$ configurations have already been classified.

There are 21,399 isomorphism classes
of $14_3$ configurations~\cite{BBP}. Checking the isomorphism classes
of graphs underlying the $14_3$ configurations
shows that their number is 20,787.
As this is about one thousandth of the number of regular graphs, the $14_3$
configurations are the appropriate building block for our algorithm.

\paragraph{Example}

There is a unique isomorphism class of an \STS{21} that contains at least
one \STS*{7} and that admits an automorphism group of order 108,
see Table~\ref{tab:result}. The following incidence matrix of such
a design visualizes the partitions of points and blocks in
the general approach (note that the ordering of rows and columns within
each subset does not necessarily coincide with the ordering given by
the algorithm):

\vspace*{5mm}
{
\tiny
\setlength{\tabcolsep}{0.65pt}
\renewcommand{\arraystretch}{1.0}
\begin{tabular}{c|ccccccc|ccccccccccccccccccccccccccccccccccccccccccccccccc|cccccccccccccc}
\multirow{7}{*}{$W$}
&1&1&1&.&.&.&.&1&1&1&1&1&1&1&.&.&.&.&.&.&.&.&.&.&.&.&.&.&.&.&.&.&.&.&.&.&.&.&.&.&.&.&.&.&.&.&.&.&.&.&.&.&.&.&.&.&.&.&.&.&.&.&.&.&.&.&.&.&.&.\\
&1&.&.&1&1&.&.&.&.&.&.&.&.&.&1&1&1&1&1&1&1&.&.&.&.&.&.&.&.&.&.&.&.&.&.&.&.&.&.&.&.&.&.&.&.&.&.&.&.&.&.&.&.&.&.&.&.&.&.&.&.&.&.&.&.&.&.&.&.&.\\
&1&.&.&.&.&1&1&.&.&.&.&.&.&.&.&.&.&.&.&.&.&1&1&1&1&1&1&1&.&.&.&.&.&.&.&.&.&.&.&.&.&.&.&.&.&.&.&.&.&.&.&.&.&.&.&.&.&.&.&.&.&.&.&.&.&.&.&.&.&.\\
&.&1&.&1&.&1&.&.&.&.&.&.&.&.&.&.&.&.&.&.&.&.&.&.&.&.&.&.&1&1&1&1&1&1&1&.&.&.&.&.&.&.&.&.&.&.&.&.&.&.&.&.&.&.&.&.&.&.&.&.&.&.&.&.&.&.&.&.&.&.\\
&.&.&1&.&1&1&.&.&.&.&.&.&.&.&.&.&.&.&.&.&.&.&.&.&.&.&.&.&.&.&.&.&.&.&.&1&1&1&1&1&1&1&.&.&.&.&.&.&.&.&.&.&.&.&.&.&.&.&.&.&.&.&.&.&.&.&.&.&.&.\\
&.&.&1&1&.&.&1&.&.&.&.&.&.&.&.&.&.&.&.&.&.&.&.&.&.&.&.&.&.&.&.&.&.&.&.&.&.&.&.&.&.&.&1&1&1&1&1&1&1&.&.&.&.&.&.&.&.&.&.&.&.&.&.&.&.&.&.&.&.&.\\
&.&1&.&.&1&.&1&.&.&.&.&.&.&.&.&.&.&.&.&.&.&.&.&.&.&.&.&.&.&.&.&.&.&.&.&.&.&.&.&.&.&.&.&.&.&.&.&.&.&1&1&1&1&1&1&1&.&.&.&.&.&.&.&.&.&.&.&.&.&.\\
\hline
\multirow{14}{*}{$V \setminus W$}
&.&.&.&.&.&.&.&1&.&.&.&.&.&.&1&.&.&.&.&.&.&1&.&.&.&.&.&.&1&.&.&.&.&.&.&1&.&.&.&.&.&.&1&.&.&.&.&.&.&1&.&.&.&.&.&.&1&1&1&.&.&.&.&.&.&.&.&.&.&.\\
&.&.&.&.&.&.&.&.&1&.&.&.&.&.&.&1&.&.&.&.&.&1&.&.&.&.&.&.&.&1&.&.&.&.&.&.&1&.&.&.&.&.&.&1&.&.&.&.&.&.&1&.&.&.&.&.&.&.&.&1&1&1&.&.&.&.&.&.&.&.\\
&.&.&.&.&.&.&.&.&.&1&.&.&.&.&.&.&1&.&.&.&.&.&1&.&.&.&.&.&.&.&1&.&.&.&.&.&.&1&.&.&.&.&.&1&.&.&.&.&.&.&.&1&.&.&.&.&1&.&.&.&.&.&1&1&.&.&.&.&.&.\\
&.&.&.&.&.&.&.&.&.&.&1&.&.&.&.&.&1&.&.&.&.&.&.&1&.&.&.&.&.&1&.&.&.&.&.&.&.&.&1&.&.&.&.&.&1&.&.&.&.&.&.&.&1&.&.&.&.&1&.&.&.&.&.&.&1&1&.&.&.&.\\
&.&.&.&.&.&.&.&.&1&.&.&.&.&.&1&.&.&.&.&.&.&.&.&.&1&.&.&.&.&.&.&1&.&.&.&.&.&.&.&1&.&.&.&.&.&1&.&.&.&.&.&.&.&1&.&.&.&.&.&.&.&.&1&.&1&.&1&.&.&.\\
&.&.&.&.&.&.&.&.&.&.&.&1&.&.&.&.&.&1&.&.&.&.&1&.&.&.&.&.&.&.&.&.&1&.&.&.&.&.&.&.&1&.&1&.&.&.&.&.&.&.&.&.&.&.&1&.&.&.&.&1&.&.&.&.&.&1&1&.&.&.\\
&.&.&.&.&.&.&.&.&.&.&.&.&1&.&.&.&.&.&1&.&.&.&.&.&.&1&.&.&.&.&.&.&.&1&.&.&.&.&.&.&.&1&.&.&.&.&1&.&.&.&.&.&.&.&.&1&1&.&.&1&.&.&.&.&1&.&.&.&.&.\\
&.&.&.&.&.&.&.&1&.&.&.&.&.&.&.&1&.&.&.&.&.&.&.&.&1&.&.&.&.&.&1&.&.&.&.&.&.&.&.&.&1&.&.&.&1&.&.&.&.&.&.&.&.&.&.&1&.&.&.&.&.&.&.&.&.&.&.&1&1&1\\
&.&.&.&.&.&.&.&.&.&.&.&.&.&1&.&.&.&.&.&1&.&.&.&.&.&1&.&.&.&.&.&.&1&.&.&.&.&1&.&.&.&.&.&.&.&1&.&.&.&.&.&.&1&.&.&.&.&.&1&.&1&.&.&.&.&.&.&1&.&.\\
&.&.&.&.&.&.&.&.&.&.&.&1&.&.&.&.&.&.&.&.&1&.&.&1&.&.&.&.&.&.&.&.&.&.&1&1&.&.&.&.&.&.&.&.&.&.&1&.&.&.&.&.&.&1&.&.&.&.&.&.&.&1&.&1&.&.&.&1&.&.\\
&.&.&.&.&.&.&.&.&.&.&1&.&.&.&.&.&.&1&.&.&.&.&.&.&.&.&1&.&.&.&.&1&.&.&.&.&.&.&.&.&.&1&.&.&.&.&.&1&.&.&1&.&.&.&.&.&.&.&1&.&.&.&.&1&.&.&.&.&1&.\\
&.&.&.&.&.&.&.&.&.&.&.&.&.&1&.&.&.&.&1&.&.&.&.&.&.&.&.&1&.&.&.&.&.&.&1&.&1&.&.&.&.&.&.&.&.&.&.&.&1&.&.&.&.&.&1&.&.&1&.&.&.&.&1&.&.&.&.&.&1&.\\
&.&.&.&.&.&.&.&.&.&1&.&.&.&.&.&.&.&.&.&.&1&.&.&.&.&.&1&.&.&.&.&.&.&1&.&.&.&.&.&1&.&.&.&.&.&.&.&.&1&1&.&.&.&.&.&.&.&.&.&.&1&.&.&.&.&1&.&.&.&1\\
&.&.&.&.&.&.&.&.&.&.&.&.&1&.&.&.&.&.&.&1&.&.&.&.&.&.&.&1&1&.&.&.&.&.&.&.&.&.&1&.&.&.&.&.&.&.&.&1&.&.&.&1&.&.&.&.&.&.&.&.&.&1&.&.&.&.&1&.&.&1\\
\hline
&&&&&&&&
\multicolumn{7}{c|}{$\mathcal{B}_0$}&
\multicolumn{7}{c|}{$\mathcal{B}_1$}&
\multicolumn{7}{c|}{$\mathcal{B}_2$}&
\multicolumn{7}{c|}{$\mathcal{B}_3$}&
\multicolumn{7}{c|}{$\mathcal{B}_4$}&
\multicolumn{7}{c|}{$\mathcal{B}_5$}&
\multicolumn{7}{c|}{$\mathcal{B}_6$}\\
\cline{9-57}
&
\multicolumn{7}{c|}{$\mathcal{B}'$}&
\multicolumn{49}{c|}{$\mathcal{F}$}&
\multicolumn{14}{c}{$\mathcal{D}$}\\
\end{tabular}
}

\vspace*{5mm}
  Let $V = \{0,1,\ldots ,20\}$. The design can be
constructed by considering the group of order 108 generated by
\begin{align*}
  &(0, 9,19)(2,10,16)(3,4,20,8,7,18)(5,6,15,14,11,13) \text{\ and}\\
  &(0,3,4)(2,5,6)(7,8, 9,20,18,19)(10,15,13,16,11,14)(12,17)
\end{align*}
and taking the 7 orbits under the action of this group
with representatives
\begin{align*}
  \{0,1,2\},\{0,3,6\},\{0, 9,19\},\{0,10,17\},\{1,12,17\},\{2,5,6\},\{2,10,16\}.
\end{align*}

\subsection{Details of the Approach}\label{sect:detail}

We shall now give more specific details needed for implementing
the general approach. Some of the computational subproblems
will be considered separately in Section~\ref{sect:sub}.

\paragraph{The point set}
When building up an \STS{21} $(V,\mathcal{B})$ containing a \STS*{7}, we let
$V = \{0,1,\ldots,20\}$ such that $W = \{14,15,\ldots,20\}$ is the point
set of the particularized \STS*{7} $(W,\mathcal{B}')$, called $S'$.

\paragraph{The $14_3$ configuration}
The distribution of the orders of the automorphism groups of the 21,399
$14_3$ configurations~\cite{BBP} is
\begin{align*}
1^{20{,}328} 2^{916} 3^{19} 4^{91} 6^{12} 7^{1} 8^{15} 12^{7} 14^{3} 16^{3} 24^{2} 128^{1} 56448^{1}.
\end{align*}
The unique $14_3$ configuration with automorphism group order 56448
consists of two disjoint \STS{7} and is the configuration of
Wilson-type systems. Ignoring that configuration here,
the groups to be considered have order at most 128,
so there is no need for advanced group algorithms.

After fixing a configuration $(V \setminus W, \mathcal D)$, where the
point set is $V \setminus W = \{0,\ldots,13\}$, we compute its automorphism
group $A$, the underlying graph $\overline{G}$, and the complement $G$.
Notice that the group $A$ is trivial in most of the cases.

\paragraph{The \mbox{1-factorization}}
For a given graph $G$, we first determine the set $F$ of \mbox{1-factors} of $G$
and then use the \mbox{1-factors} in $F$ to compute the set $\mathcal{F}'$ of
all possible \mbox{1-factorizations} of $G$.
If the group $A$ is nontrivial, isomorph rejection is further carried
out by accepting precisely those \mbox{1-factorizations} in $\mathcal{F}'$ that are
lexicographically minimum under the action of $A$. For an
accepted \mbox{1-factorization}, the subgroup of
$A$ consisting of the elements that
stabilize the \mbox{1-factorization} is denoted by $A'$.

A \mbox{1-factor} of $G$ corresponds to a set $\mathcal{B}'_p$
defined in~\eqref{eq:1}, and a \mbox{1-factorization} of $G$ gives a set of
blocks $\mathcal{F} = \cup_{p=14}^{20}\mathcal{B}_p$ up to permutation of the
points in $W = \{14,15,\ldots ,20\}$ (we pick an arbitrary one).
The group $A'$ acts on $V \setminus W$. Blocks of $\mathcal{F}$
also have points in $W$, so we extend the action of $A'$ to get a
group $A''$ acting on $V$. The permutation of the points in
$W$ for an element in $A''$ is uniquely defined by how the original
element in $A'$ maps the \mbox{1-factors}.

\paragraph{The sub-STS(7)}

There is a unique \STS{7}, the Fano plane, which has an automorphism
group of order 168. Hence there are are $7!/168 = 30$ distinct
labelled $\STSS{7}$ on 7 given points.

An isomorphism from one \STS{21} with a \STS*{7} to another
maps the particularized \STS*{7} to a \STS*{7}. Hence there are two
general situations: \STSS{21} with exactly one \STS*{7} and
\STSS{21} with more than one \STSS*{7}. In the latter case,
there are further several possibilities for how the point
sets of two \STSS*{7} may intersect. Such an intersection must form
a (possibly trivial) sub-STS, so possible intersection sizes are
0, 1, and 3.

If the intersection size is 0, then there is necessarily a
third \STS*{7} whose point set is disjoint from the point sets
of the first two \STSS*{7}, that is, we have a Wilson-type system
and the $14_3$ configuration discussed earlier. Wilson-type
\STSS{21} have exactly three \STSS*{7}~\cite[Lemma~1]{KOTZ}.
As the mentioned $14_3$ configuration is not considered here, this case will
not occur in the search.

Isomorph rejection when extending blocks $\mathcal{D} \cup \mathcal{F}$
with blocks $\mathcal{B}'$ is analogous to the situation when extending
blocks $\mathcal{D}$ with blocks $\mathcal{F}$, considered earlier. Now, out of
the 30 possibilities, those \STSS*{7} that are lexicographically
minimum under the action of $A''$ are accepted.
The subgroup of $A''$ consisting of the elements that
stabilize the accepted \STS*{7} is denoted by $A'''$.

The blocks $\mathcal{B} = \mathcal{D} \cup \mathcal{F} \cup
\mathcal{B}'$ now form an \STS{21} with a \emph{particularized} \STS*{7},
and if those are the objects to classify we would be done. But in
the classification of \STSS{21} with \emph{at least one} \STS*{7},
there is still one final step.

\paragraph{The final isomorph rejection}

If there is exactly one \STS*{7} in the constructed design $(V,\mathcal{B})$,
then we accept the \STS{21}; its automorphism group is the group $A'''$
computed earlier. Otherwise, we proceed by finding all \STSS*{7} in $V$.
(As we have seen, these will intersect $W$ in exactly
1 or 3 points;
some precomputations for finding them can be done based on $\mathcal{D}$
and $\mathcal{F}$.)
We now determine whether the particularized \STS*{7} is a canonically
minimum \STS*{7}, to be discussed in Section~\ref{sect:sub}, and accept it
if that is the case. The automorphism group of an accepted \STS{21} is
obtained as a by-product of the computations.

\subsection{Computational Subproblems}\label{sect:sub}

We shall here discuss some of the main computational subproblems that are
encountered when implementing the presented approach and that are not
standard problems related to data structures and algorithms.

\paragraph{Automorphism groups and canonical forms}
Automorphism groups and canonical forms are conveniently computed
with \texttt{nauty}~\cite{MP} after an appropriate transformation
of the combinatorial structure to a graph.

To order the \STSS*{7} of an \STS{21} one may
use the standard graph encoding of the incidence matrix of the
design, add one vertex for each \STS*{7}, and let the 7 vertices
corresponding to the points of the \STS*{7} form the neighborhood
of an added vertex. Then the canonical order of vertices given by
\texttt{nauty} imposes an order on the \STSS*{7}. More precisely,
\texttt{nauty} determines an order of the \emph{orbits} of vertices under the
action of the automorphism group of the graph. Therefore we get an induced ordering of
the orbits of \STSS*{7} under the action of the automorphism group of the
\STS{21}.

For small group orders, the abstract type of the automorphism groups
of the classified designs can be identified
based on the multiset of orders of elements.
The abstract type can further be computed using
\texttt{AllSmallGroups} and \texttt{StructureDescription} in
\texttt{GAP}~\cite{GAP}. In the current work, seven groups (of orders
27, 54, 108, 294, and 1008) had to be treated manually and separately.
The designs with nontrivial automorphisms are amongst those
classified in~\cite{K1}.

\paragraph{\mbox{1-factors} and \mbox{1-factorizations}}
Although finding single perfect matchings in general graphs respective
all perfect matchings in bipartite graphs are standard computational
problems~\cite{E,FM}, we use a backtrack algorithm to compute all
\mbox{1-factors} of general graphs.
Given the set of \mbox{1-factors} of a
graph, the problem of finding all \mbox{1-factorizations} can be phrased in the
framework of exact cover~\cite{KO3}, whereby the instances can be solved,
for example, using the \texttt{libexact}~\cite{KP} software.

\section{Results}\label{sect:res}

The total number of isomorphism classes of
\STSS{21} containing at least one \STS*{7} is 116,635,963,205,551,
which splits into 116,635,961,039,200 cases that are not
of Wilson type and 2,166,351 cases that are of Wilson type~\cite{KOTZ}.

More detailed information can be found in Table~\ref{tab:result}
and Table~\ref{tab:aggregated}.
The column headers in Table~\ref{tab:result} are the order of the
automorphism group ($O$), the number of contained
\STSS*{7} ($U$), the number of unordered pairs of \STSS*{7} that
intersect in 1 and 3 points ($I_{1}$ and $I_{3}$, respectively),
the abstract type of the automorphism group ($A$), and finally the
number of isomorphism classes of \STSS{21} with these properties ($\#$).

For completeness, Table~\ref{tab:W} from~\cite{KOTZ} is included.
For all Wilson-type \STSS{21}, we have $U=3$, $I_{1}=0$, and
$I_{3}=0$ by~\cite[Lemma~1]{KOTZ}.
In the Appendix, some data for \STSS{21} that do not contain \STSS*{7} is given.

The notation for the abstract types of groups is as follows:
$C_{n}$ is the cyclic group of order $n$,
$S_{n}$ is the symmetric group of order $n!$,
$A_{n}$ is the alternating group of order $n!/2$,
$D_{n}$ is the dihedral group of order $n$, and
$\operatorname{PSL}(v,q)$ is the projective special linear group in
$\mathbb{F}_q^v$.
For two groups $G$ and $H$, $G \times H$ is the direct product of
$G$ and $H$, $G \rtimes H$ is a semidirect product of $G$ and $H$, and
$G^n$ is $G \times G \times \cdots \times G$ ($n$ times).

A central open problem for specific \STSS{21} is whether systems
exist that are doubly resolvable.
The current work gives nothing
new with respect to this problem, because Kirkman triple systems
of order 21 with
\STSS*{7} have already been classified and tested~\cite{KO5}.

The whole classification including the detection of the abstract group
types took
about 1300 CPU days on the equivalent of one core of an
Intel Xeon E5-2665 @ 2.40GHz.

\paragraph{Verification}
We perform two tests to validate results.
Let $\mathcal{S}$ be a transversal of the isomorphism classes of
the \STSS{21} with
\STSS*{7} that are not of Wilson type---this is the outcome
of the current classification---and let $\mathcal{C}$ be a transversal
of the isomorphism classes of the $14_3$ configurations excluding the
configuration leading to Wilson-type \STSS{21}. Further, let
$s_7(S)$ be the number of \STSS*{7} in the system $S$,
and let $f(C)$ be the number of
\mbox{1-factorizations} with labelled \mbox{1-factors} of the complement of the graph
underlying the configuration $C$.
During the computations, all this data was collected.

In the first test, we count in two different ways all pairs of
labelled \STSS{21}  that are not of Wilson type and their contained
\STSS*{7}. By the Orbit--Stabilizer Theorem, we have

\begin{align*}
\sum_{S \in \mathcal{S}} \frac{21!}{|\operatorname{Aut}(S)|} \cdot s_7(S)
=
\sum_{C \in \mathcal{C}} \frac{21!}{|\operatorname{Aut}(C)|} \frac{7!}{168}\cdot f(C)
.
\end{align*}

Both sides of this equality yielded
\[
5{,}988{,}986{,}139{,}804{,}614{,}556{,}727{,}954{,}636{,}800{,}000
\]
in the final computation.
The fact~\cite{KOTZ} that Wilson-type \STSS{21} have
exactly three \STSS*{7} and will not appear in the search is essential
for the double counting to work.

In the second test, we extract the \STSS{21} with \STSS*{7}
from the \STSS{21} with nontrivial automorphisms classified in~\cite{K1}
and compare the numbers with those in Table~\ref{tab:nonW}. Also
this test was successful.

\begin{table}
\caption{Numbers of \protect\STSS{21} containing at least one \protect\STS*{7}\label{tab:result}}
\centering
\scriptsize
\setlength{\tabcolsep}{2pt}
\renewcommand{\arraystretch}{1.5}
\begin{subtable}[t]{0.7\textwidth}
\caption{non-Wilson type\label{tab:nonW}}
\begin{minipage}{0.52\textwidth}
\centering
\begin{tabular}{rrrrlr}
\toprule
$O$ & $U$ & $I_{1}$ & $I_{3}$ & $A$ & $\#$ \\
\midrule
1 & 1 & 0 & 0   & $C_{1}$             & 116,051,875,827,936 \\
1 & 2 & 1 & 0   & $C_{1}$             & 31,778,146,776 \\
1 & 2 & 0 & 1   & $C_{1}$             & 550,238,290,596 \\
1 & 3 & 1 & 2   & $C_{1}$             & 593,663,600 \\
1 & 3 & 3 & 0   & $C_{1}$             & 60,352,088 \\
1 & 3 & 0 & 3   & $C_{1}$             & 1,385,739,943 \\
1 & 4 & 1 & 5   & $C_{1}$             & 6,391,040 \\
1 & 4 & 2 & 4   & $C_{1}$             & 198,304 \\
1 & 4 & 3 & 3   & $C_{1}$             & 1,607,028 \\
1 & 4 & 0 & 6   & $C_{1}$             & 157,886 \\
1 & 5 & 1 & 9   & $C_{1}$             & 576 \\
1 & 5 & 2 & 8   & $C_{1}$             & 50,192 \\
1 & 5 & 3 & 7   & $C_{1}$             & 30,024 \\
1 & 5 & 4 & 6   & $C_{1}$             & 1,704 \\
1 & 6 & 3 & 12  & $C_{1}$             & 1,790 \\
1 & 6 & 4 & 11  & $C_{1}$             & 688 \\
1 & 7 & 5 & 16  & $C_{1}$             & 124 \\
2 & 1 & 0 & 0   & $C_{2}$             & 19,270,679 \\
2 & 2 & 1 & 0   & $C_{2}$             & 84,080 \\
2 & 2 & 0 & 1   & $C_{2}$             & 814,880 \\
2 & 3 & 1 & 2   & $C_{2}$             & 18,912 \\
2 & 3 & 3 & 0   & $C_{2}$             & 43,062 \\
2 & 3 & 0 & 3   & $C_{2}$             & 132,334 \\
2 & 4 & 1 & 5   & $C_{2}$             & 9,088 \\
2 & 4 & 2 & 4   & $C_{2}$             & 64 \\
2 & 4 & 3 & 3   & $C_{2}$             & 2,448 \\
2 & 5 & 2 & 8   & $C_{2}$             & 224 \\
2 & 5 & 3 & 7   & $C_{2}$             & 2,092 \\
2 & 5 & 4 & 6   & $C_{2}$             & 16 \\
2 & 6 & 3 & 12  & $C_{2}$             & 140 \\
2 & 6 & 4 & 11  & $C_{2}$             & 32 \\
2 & 7 & 5 & 16  & $C_{2}$             & 188 \\
2 & 9 & 9 & 27  & $C_{2}$             & 2 \\
3 & 1 & 0 & 0   & $C_{3}$             & 177,205 \\
3 & 2 & 1 & 0   & $C_{3}$             & 3,152 \\
3 & 2 & 0 & 1   & $C_{3}$             & 5,508 \\
\bottomrule
\end{tabular}
\end{minipage}
\begin{minipage}{0.5\textwidth}
\begin{tabular}{rrrrlr}
\toprule
\hspace*{2mm}$O$ & $U$ & $I_{1}$ & $I_{3}$ & $A$ & $\#$\\
\midrule
3 & 3 & 3 & 0   & $C_{3}$             & 655 \\
3 & 3 & 0 & 3   & $C_{3}$             & 4,152 \\
3 & 4 & 3 & 3   & $C_{3}$             & 132 \\
3 & 4 & 0 & 6   & $C_{3}$             & 6 \\
3 & 5 & 4 & 6   & $C_{3}$             & 16 \\
3 & 6 & 3 & 12  & $C_{3}$             & 18 \\
4 & 1 & 0 & 0   & $C_{2}^2$             & 6,268 \\
4 & 1 & 0 & 0   & $C_{4}$             & 628 \\
4 & 3 & 3 & 0   & $C_{2}^2$             & 260 \\
4 & 3 & 0 & 3   & $C_{2}^2$             & 870 \\
4 & 5 & 3 & 7   & $C_{2}^2$             & 136 \\
4 & 7 & 5 & 16  & $C_{2}^2$             & 24 \\
4 & 9 & 9 & 27  & $C_{2}^2$             & 3 \\
6 & 1 & 0 & 0   & $C_{6}$             & 849 \\
6 & 1 & 0 & 0   & $S_{3}$             & 192 \\
6 & 3 & 3 & 0   & $C_{6}$             & 146 \\
6 & 3 & 3 & 0   & $S_{3}$             & 39 \\
6 & 3 & 0 & 3   & $C_{6}$             & 91 \\
6 & 3 & 0 & 3   & $S_{3}$             & 31 \\
6 & 4 & 3 & 3   & $S_{3}$             & 16 \\
6 & 6 & 3 & 12  & $S_{3}$             & 4 \\
6 & 9 & 9 & 27  & $S_{3}$             & 2 \\
7 & 1 & 0 & 0   & $C_{7}$             & 27 \\
8 & 1 & 0 & 0   & $C_{4} \times C_{2}$         & 8 \\
8 & 1 & 0 & 0   & $D_{8}$             & 164 \\
9 & 3 & 3 & 0   & $C_{3}^2$             & 1 \\
9 & 3 & 0 & 3   & $C_{3}^2$             & 3 \\
12  & 3 & 3 & 0   & $D_{12}$              & 4 \\
12  & 3 & 0 & 3   & $D_{12}$              & 10 \\
12  & 9 & 9 & 27  & $D_{12}$              & 3 \\
14  & 1 & 0 & 0   & $C_{14}$              & 14 \\
16  & 1 & 0 & 0   & $C_{2} \times D_{8}$         & 8 \\
18  & 3 & 3 & 0   & $C_{3} \times S_{3}$         & 11 \\
18  & 3 & 0 & 3   & $C_{3} \times S_{3}$         & 6 \\
36  & 9 & 9 & 27  & $S_{3}^2$         & 1 \\
108 & 9 & 9 & 27  & $(C_{3}^2 \rtimes C_{6}) \rtimes C_{2}$   & 1 \\
\bottomrule
\end{tabular}
\end{minipage}
\end{subtable}
\hspace*{2.3mm}
\begin{subtable}[t]{0.2\textwidth}
\caption{Wilson type\label{tab:W}}
\centering
\begin{tabular}{rlr}
\toprule
\hspace*{6mm} $O$ & $A$ & $\#$ \\
\midrule
1 		& $C_{1}$ & 2,156,186 \\
2 		& $C_{2}$ & 8,914 \\
3 		& $C_{3}$ & 685 \\
4 		& $C_{2}^2$ & 253 \\
4 		& $C_{4}$ & 18 \\
6 		& $C_{6}$ & 94 \\
6 		& $S_{3}$ & 103 \\
8 		& $C_{2}^3$ & 22 \\
8 		& $D_{8}$ & 19 \\
9 		& $C_{3}^2$ & 3 \\
12 		& $A_{4}$ & 2 \\
12 		& $D_{12}$ & 4 \\
16 		& $C_{2} \times D_{8}$ & 4 \\
18 		& $C_{3}^2 \rtimes C_{2}$ & 2 \\
18 		& $C_{3} \times S_{3}$ & 5 \\
21 		& $C_{7} \rtimes C_{3}$ & 2 \\
24 		& $C_{2} \times A_{4}$ & 9 \\
24 		& $C_{2}^2 \times S_{3}$ & 1 \\
24 		& $S_{4}$ & 7 \\
42 		& $C_{2} \times (C_{7} \rtimes C_{3})$ & 1 \\
42 		& $C_{7} \rtimes C_{6}$ & 5 \\
48 		& $C_{2} \times S_{4}$ & 2 \\
72 		& $(C_{3} \times A_{4}) \rtimes C_{2}$ & 1 \\
72 		& $A_{4} \times S_{3}$ & 4 \\
126 	& $S_{3} \times (C_{7} \rtimes C_{3})$ & 1 \\
144 	& $S_{3} \times S_{4}$ & 1 \\
294 	& $C_{7}^2 \rtimes C_{6}$ & 1 \\
882 	& $(C_{7} \rtimes C_{3})^2 \rtimes C_{2}$ & 1 \\
1008 	& $\operatorname{PSL}(3,2) \times S_{3}$ & 1 \\
\bottomrule
\end{tabular}
\end{subtable}
\end{table}

\begin{table}
\caption{Aggregated numbers of \protect\STSS{21} containing at least one \protect\STS*{7}\label{tab:aggregated}}
\hfill
\begin{tabular}{rr}
\toprule
$O$ & $\#$ \\
\midrule
1 & 116,635,942,616,481 \\
2 & 20,387,155 \\
3 & 191,529 \\
4 & 8,460 \\
6 & 1,567 \\
7 & 27 \\
\bottomrule
\end{tabular}
\hfill
\begin{tabular}{rr}
\toprule
$O$ & $\#$ \\
\midrule
8 & 213 \\
9 & 7 \\
12 & 23 \\
14 & 14 \\
16 & 12 \\
18 & 24 \\
\bottomrule
\end{tabular}
\hfill
\begin{tabular}{rr}
\toprule
$O$ & $\#$ \\
\midrule
21 & 2 \\
24 & 17 \\
36 & 1 \\
42 & 6 \\
48 & 2 \\
72 & 5 \\
\bottomrule
\end{tabular}
\hfill
\begin{tabular}{rr}
\toprule
$O$ & $\#$ \\
\midrule
108 & 1 \\
126 & 1 \\
144 & 1 \\
294 & 1 \\
882 & 1 \\
1008 & 1 \\
\bottomrule
\end{tabular}
\hfill
\end{table}

\section{Estimating the Number of STS(21)s}\label{sect:estimate}

The classification of the \STSS{21} with \STSS*{7} gives a lower bound on
the number of isomorphism classes of \STS{21} but can
also be used for an estimation of that number.
The authors are not aware of any published estimations.

Quackenbush~\cite{Q} conjectured that almost all Steiner triple systems have no nontrivial subsystems.
Later, however, Kwan~\cite{K0} used a random model to find evidence for the number of \STSS*{7} in an \STS{v} to have expectation $\Theta(1)$, referring to similar work in~\cite{MW} on Latin squares.
The models used in~\cite{K0} and~\cite{MW} are random \mbox{3-uniform} hypergraphs
and random integer matrices, respectively.

An \STS{v} has $v(v-1)/6$ blocks out of $v(v-1)(v-2)/6$ \mbox{3-subsets} of a $v$-set, that is, a ratio of $p := 1/(v-2)$ of the \mbox{3-subsets} are blocks.
We may now form a random \mbox{3-uniform} hypergraph on $v$ vertices by including blocks with probability $p$ (note that $p := 1/v$, which is used in~\cite{K0}, works when studying asymptotics). We denote the
number of labelled \STSS{w} on $w$ points by $N(w)$.
We have seen earlier that $N(7)=30$.
The number of labelled \STSS{w} on $v$ points, where $v \geq w$, is
$M(v,w):=N(w)\binom{v}{w}$. The probability for a given \STS{w} to occur
in the random model is $p^{w(w-1)/6}$.

The linearity of the expected value allows now to compute the expected number of \STS*{w}
\begin{align*}
\mu(v,w) := \frac{N(w)\binom{v}{w}}{(v-2)^{w(w-1)/6}}
\end{align*}
and, abbreviating $\mu(\infty,w) = \lim_{v\to\infty} \mu(v,w)$, we have $\mu(\infty,7) = 1/168 \approx 0.00595$ and $\mu(\infty,w) = 0$ for $w > 7$.

Let $S$ be the set of positive integers fulfilling~\eqref{eq:sts}.
Analogously to the conjecture in~\cite[p.~346]{MW}, see also~\cite{K0}, we state
the following.

\begin{conjecture}
The distribution of the number of \STSS*{w} in an \STS{v} tends to the
Poisson distribution with expected value $1/168$ for $w=7$ and expected
value $0$ for $w>7$ as $v \in S$ tends to infinity.
\end{conjecture}

The proportion of \STS{v} containing at least one \STS*{7} is then approximately $\alpha = 1-e^{-1/168} \approx 0.00593$ for large $v$.
Consequently, an estimation of the total number of \STS{v} can be obtained by dividing the number of \STS{v} with at least one \STS*{w} by $\alpha$.

As almost all Steiner triple systems have no nontrivial automorphisms~\cite{B}, an estimation for the number of isomorphism classes of \STS{v} can be obtained by dividing the number of isomorphism classes of \STS{v} with at least one \STS*{w} by $\alpha$.

There are only two instances for which the quality of such an estimation can be checked, and the case \STS{15} with only 80 isomorphism classes is not useful.
For the 11,084,874,829 \STS{19} out of which 86,701,547 have at least one \STS*{7}~\cite{KOTZ}, we get a ratio of approximately 0.00782, which is getting into the same order of magnitude as $\alpha$.

For the number of isomorphism classes of \STSS{21}, using
the classification results of the current paper we calculate

\[
116{,}635{,}963{,}205{,}551 / \alpha \approx 1.965 \cdot 10^{16},
\]

\noindent
which indicates that the number could be somewhat greater than
$10^{16}$, perhaps between $1\cdot 10^{16}$ and $2\cdot 10^{16}$.

In the estimation one might consider utilizing $\mu(21,7) \approx 0.00389$
rather than $\mu(\infty ,7)$, but notice that
$\mu(19,7) \approx 0.00368$ underestimates the true value by a factor
greater than 2,
and the situation here might be analogous to that for sub-Latin squares
considered in~\cite{MW}.
In that paper, it is conjectured that the expected number of $3 \times 3$ sub-Latin squares of a randomly chosen
$n \times n$ Latin square tends to $1/18$ as $n$ tends to infinity,
and numerical data show that the value given by the random model
for a fixed parameter, $f(n)=12 \binom{n}{3}^3 n^{-9}$, underestimates the computed value for small parameters. For example, for $n=10$, the asymptotic value ($\approx 0.0556$) is closer to the exact value ($\approx 0.0536$) than $f(n)$ ($\approx 0.0207$).

It is not clear whether an STS with subsystems is more or less
prone to have resolutions. If the correlation is weak, then the
fact that there are 12,520,021 isomorphism classes of Kirkman
triple systems of order 21 with $\STSS*{7}$~\cite{KO5}
could be used to calculate
\[
12{,}520{,}021 / \alpha \approx 2.111 \cdot 10^{9},
\]
which would hint that there might be somewhat more than
1 billion isomorphism classes of Kirkman triple systems of order 21.

\section*{Appendix}

The program developed for determining the abstract type of an
automorphism group can be applied also to those \STSS{21} with
nontrivial automorphisms from~\cite{K1} that do not contain \STSS*{7}.
Such information is presented in Table~\ref{tab:other}
using the notation described in Section~\ref{sect:res}.

\begin{table}
\caption{Numbers of \protect\STSS{21} with no \protect\STSS*{7}}\label{tab:other}
\hfill
\begin{tabular}{rlr}
\toprule
$O$ & $A$ & $\#$ \\
\midrule
2 		& $C_{2}$ 									& 40,201,112 \\
3 		& $C_{3}$ 									& 1,540,602 \\
4 		& $C_{2}^2$ 								        & 3,007 \\
5 		& $C_{5}$ 									& 1,772 \\
6 		& $C_{6}$ 									& 533 \\
6 		& $S_{3}$ 									& 279 \\
7 		& $C_{7}$ 									& 39 \\
8 		& $D_{8}$ 									& 9 \\
9 		& $C_{3}^2$ 								& 95 \\
9 		& $C_{9}$ 									& 7 \\
12 		& $A_{4}$ 									& 44 \\
12 		& $C_{6} \times C_{2}$ 						& 18 \\
18 		& $C_{3}^2 \rtimes C_{2}$ 					& 1 \\
\bottomrule
\end{tabular}
\hfill
\begin{tabular}{rlr}
\toprule
$O$ & $A$ & $\#$ \\
\midrule
18 		& $C_{18}$ 									& 2 \\
18 		& $C_{3} \times S_{3}$ 						& 1 \\
18 		& $C_{6} \times C_{3}$ 						& 5 \\
21 		& $C_{21}$ 									& 1 \\
21 		& $C_{7} \rtimes C_{3}$ 					& 7 \\
24 		& $C_{3} \times D_{8}$ 						& 1 \\
24 		& $S_{4}$ 									& 1 \\
27 		& $C_{3}^2 \rtimes C_{3}$ 					& 3 \\
36 		& $C_{3} \times A_{4}$ 						& 4 \\
42 		& $C_{3} \times D_{14}$ 					& 1 \\
54 		& $C_{3}^2 \rtimes C_{6}$ 					& 1 \\
126 	& $C_{3} \times (C_{7} \rtimes C_{6})$ 		& 1 \\
504 	& $C_{3} \times \operatorname{PSL}(3,2)$ 	& 1 \\
\bottomrule
\end{tabular}
\hfill
\end{table}

\section*{Acknowledgements}

The authors are grateful to Petteri Kaski for providing the \STSS{21}
with nontrivial automorphisms classified in~\cite{K1}.

\end{document}